\documentclass{amsart}
\usepackage{graphicx}

\theoremstyle{definition}

\theoremstyle{remark}

\numberwithin{equation}{section}

\begin{document}

\title{Dynamical symmetry breaking described by cubic nonlinear Klein-Gordon equations}

\author{Yasuhiro Takei$^{1}$
}

\author{Yoritaka Iwata$^{2,*}$}

\address{$^1$Mizuho Information \& Research Institute }
\address{$^2$Osaka University of Economics and Law}
\email{y-iwata@s.keiho-u.ac.jp}





\begin{abstract}
The dynamical symmetry breaking associated with the existence and non-existence of breather solutions is studied.
Here, nonlinear hyperbolic evolution equations are calculated using a high-precision numerical scheme.
First, for clarifying the dynamical symmetry breaking, it is necessary to use a sufficiently high-precision scheme in the time-dependent framework.
Second, the error of numerical calculations is generally more easily accumulated for calculating hyperbolic equations rather than parabolic equations. 
Third, numerical calculations become easily unstable for nonlinear cases.
Our strategy for the high-precision and stable scheme is to implement the implicit Runge-Kutta method for time, and the Fourier spectral decomposition for space.
In this paper, focusing on the breather solutions, the relationship between the velocity, mass, and the amplitude of the perturbation is clarified. 
As a result, the conditions for transitioning from one state to another are clarified. 
\end{abstract}

\maketitle

\section{Background}
Nonlinearity is the main driving force of complexity (for an outline, for example, see Ref.~\cite{2011mitchell}).
Pattern formations, shock waves, solitons, and breather waves (for the typical time evolution of the breather solution, see Fig.~4 of Ref.~\cite{Y.Takei2021}, while for the time evolution of ordinary sine wave, see Fig.~3 of Ref.~\cite{Y.Takei2021}) are good examples~\cite{M.C.Cross1993,H.Nabika2020,P.Ball2015}. 
Among several nonlinear phenomena, here we focus on the transition between different states in which the existence of a breather solution is associated with the suppression of such transitions.
The existence of breather solutions was studied for several nonlinear Klein-Gordon equations; for example, Refs.~\cite{J.Denzler1993,C.Blank2011,D.Maier2018,D.Scheider2020}. 
Although breather solutions were reported to exist for the Sine-Gordon equation, they are turned out to be not necessarily stable. 
The existence of breather solutions is also known for the third-order nonlinear Klein-Gordon equation with spatial periodic boundary conditions~\cite{Y.Takei2021}. 
Note that, based on high-precision numerical calculations, the appearance of breather waves is suggested to be associated with the Lyapunov-stability constants of steady solutions~\cite{Y.Iwata2023}.

To understand the complexity, nonlinear partial differential equations have been studied in both pure mathematical and numerical ways.
For calculating nonlinear partial differential equations numerically, precision and stability must be realized simultaneously.
Furthermore, for calculating the non-stationary problems (time-dependent problems) numerically, a more rigorous criterion for the precision is generally required rather than calculating stationary problems.
Generally speaking, the simple implementation of finite differential methods for both time and space cannot lead to the correct answer for time-dependent and nonlinear hyperbolic partial differential equations.

For the third-order nonlinear Klein-Gordon equations, we employ our original high-precision numerical scheme \cite{Y.Iwata2020,Y.Takei2022} consisting of the Fourier spectral method for space and the implicit Runge-Kutta method for time. 
For applications of the high-precision numerical scheme to coupled Klein-Gordon equations, see Ref. \cite{Y.Takei2023}.
In this paper, the transitions from one Lyapunov-stable state to another stable state are time-dependently calculated for the first time in the nonlinear Klein-Gordon framework.
As a result, the appearance condition of dynamical symmetry breaking is shown that depends on the physical quantity such as the phase velocity and the initial wave amplitude.

\section{Mathematical models}
\subsection{Partial differential equations}
We consider the initial and boundary value problem (KG) of nonlinear Klein-Gordon equations with third-order nonlinear terms.
\[
\begin{array}{lll}
\displaystyle \frac{\partial^{2} u}{\partial t^{2}} + \alpha \frac{\partial^{2} u}{\partial x^{2}} + \beta u(u^{2} - \mu) = 0,
\end{array}
\hspace{5.5cm} {\rm (KG)} 
\]
where real coefficients satisfy $\alpha< 0$, and $\beta,\mu > 0$.
The new variables $\xi = x/L$ and $\tau = \beta^{1/2} t$ are introduced in a simpler and equivalent form. 
By imposing the conditions, the initial and boundary value problem (IBVP) is obtained.
\[
\begin{array}{lll}
\displaystyle \frac{\partial^{2} u}{\partial t^{2}} =  \left( \frac{-\alpha}{\beta L^{2}} \right) \frac{\partial^{2} u}{\partial x^{2}} - u(u^{2} - \mu), 
\quad x \in [0,\ 1],
\quad t \in [0,\ \infty),\ \vspace{3mm} \\
u(x, 0) = A \sin(2 \pi x) + \mu^{1/2},
\quad u(0, t) = u(1, t),  \vspace{3mm} \\
\displaystyle \frac{\partial u}{\partial x}(x, 0) = 0,
\quad \frac{\partial u}{\partial x}(0, t) = \frac{\partial u}{\partial x}(1, t).
\end{array}
\hspace{4em} {\rm (KG')} 
\]
where the variables $\xi$ and $\tau$ are still denoted by $x$ and $t$ respectively, if there is no ambiguity.
In this equation, $u=0$, $\pm \mu^{1/2}$ are steady-state constant solutions, and $u=\pm \mu^{1/2}$ are Lyapunov-stable solutions in the cases where $\mu >0$ is satisfied. 
For the initial state, a spatially-inhomogeneous perturbation $A \sin(2 \pi x)$ is added to the Lyapunov-stable state $\mu^{1/2}$, where a real constant $A$ satisfies $A \ge 0$.
As the definition of Lyapunov-stability and as shown for some examples~\cite{Y.Iwata2023}, it is expected that non-stationary solutions $u(x,t)$ keep staying around $u= \mu^{1/2}$ or $u= - \mu^{1/2}$, if $u(x,0)$ is sufficiently close to $\mu^{1/2}$ or $-\mu^{1/2}$, respectively.

\subsection{Reduced ordinary differential equations}
The reduced initial value problem (IVP) of the ordinary differential equation is derived by assuming $\alpha = 0$ to the IBVP (KG'): 
\begin{equation}  \label{ODE}
\begin{array}{lll}
\displaystyle \frac{d^{2} u}{d t^{2}} =  - u(u^{2} - \mu), 
\quad t \in [0,\ \infty),\ \vspace{3mm} \\
u(0) = A + \mu^{1/2}.  
\end{array}
\end{equation}
This reduction enables us to have a point-wise treatment of waves, although the original model (KG') takes into account the spatial distribution (finite-size effect) of waves.
That is, the comparison between the results from (KG') and (\ref{ODE}) shows the effects of spatial distribution and the resulting wave property.

\section{Results}
\subsection{Calculation of dynamical symmetry breaking}
Dynamical symmetry breaking: the transition
from one Lyapunov-stable state to another stable state is calculated systematically.
The constant values $\beta$ and $L$ are fixed to $\beta=1$ and $L=1$ respectively, while the various values of the phase velocity $\sqrt{-\alpha}$ and the wave amplitude $A$ are examined for four different mass parameters $\mu=2^{-6}$, $2^{-2}$, $2^{0}$, and $2^{1}$.
There are two possible oscillations:
\def\theenumi{\roman{enumi}}
\def\labelenumi{\theenumi)}

\begin{enumerate}
    \item oscillations confined only within the positive or negative side, \\
    \item oscillations not confined within positive or negative side.
\end{enumerate}
In the present initial values, the former cases correspond to oscillations that keep staying around the stationary solution $+ \mu^{1/2}$, and the latter cases to oscillations that do not keep staying around $+ \mu^{1/2}$.
The simple explanation and the links to the movies are found in Ref.~\cite{Y.Takei-Web}.

\begin{figure}[bt]  
\includegraphics{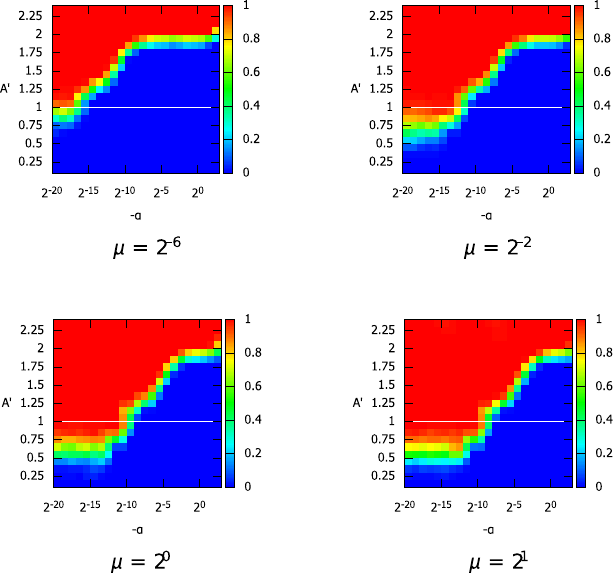} 
  \caption{\label{fig:011} 
(Color online) Phase diagram of stable oscillation based on IBVP (KG') with several $\mu$.
The horizontal axis represents $-\alpha$ in the logarithmic scale, and the vertical axis does the normalized value $A' = A /\{ (\sqrt{2}-1)  \mu^{1/2} \}$ of the original $A$. 
Based on IVP (\ref{ODE}), $A'=1$ line (white line in each panel) means the boundary value between the final states i) and ii).
Each diagram consists of 19$\times$24 pixels (squares), and each pixel is colored based on 32 stochastically chosen calculation data.
The representative values are calculated by "the number of case ii) divided by 32".
}
\end{figure}


$19$ different cases are taken from the interval $[2^{-20},2^{-2}]$ for the value of $-\alpha$. 
The interval $[0.04,0.13]$ for the value of $A$ is divided into 24 parts, and 32 values of $A$ are randomly chosen in each part. 
Consequently, systematic numerical simulations are carried out with cases
$19 \times 24 \times 32 = 14592$ to ensure statistical sufficiency.
All calculations are performed from $t = 0$ to $16384$, which is sufficient to identify the type i) or type ii) of oscillations.
The comparable calculations are performed by IVP (\ref{ODE}), which corresponds to point-wise treatment of the dynamics with spatially homogeneous perturbation.

Based on the identification of the final states, a value "0" is assigned for cases i), and a value "1" is assigned for all the other cases ii). 
Figure \ref{fig:011} summarizes the oscillation statistics: ``the number
of case ii) divided by the number of all cases”.
Since Fig. \ref{fig:011} is based on time-dependent calculations with perturbed initial states, it actually shows the dynamical breaking of the symmetry.
The red area indicates the ordinary case ii), while the blue area indicates the cases dominated by the oscillations staying around the stationary solution $+ \mu^{1/2}$.
In the intermediate yellow and green areas, the coexistence of two different types of oscillations appears.

According to Fig. \ref{fig:011}, 
oscillations confined only to the positive or negative side tend to be found in smaller amplitude waves ($A \ll 1$) and waves with larger phase velocity $\sqrt{-\alpha}$.
It can also be seen that the boundary value of $A$ does not change significantly even when $-\alpha > 2^{-6}$. 
Since the $x$-axis is a log scale, the boundary looks like a step function.
Among others, it is remarkable that the coexistence of two different oscillations is calculated.
It is , so to say, a chaotic area in which the drastically different final states are obtained due to the tiny difference of the initial conditions (cf. the sensitivity to the initial conditions).
The coexistence of two different oscillations tends to be found in smaller amplitude waves ($A \ll 1$) and waves with smaller phase velocity $\sqrt{-\alpha}$.

\section{Summary}
The dynamical symmetry breaking is shown by the existence and non-existence of a transition from a stable state to another.
In other words, it is the sustainability of the pure states, in contrast to the formation of the mixed states.
In order to realize the transition from one state to another, it has been clarified that sufficient amplitude of perturbation is required.
This tendency is also clarified to be dependent on the amplitude of phase-velocity.
According to the comparison between the results of IBVP (KG') and IVP (\ref{ODE}), the enhancement of the transition is calculated for a small phase velocity $\sqrt{-\alpha}$, and the reduction of the transition is calculated for a larger phase velocity.
This enhancement almost vanishes for a very small $\mu$ that satisfies $\mu < 2^{-6}$.

In conclusion, based on a systematic calculation, it is confirmed that the transition between the two stable states is suppressed for smaller $\mu$ ($\sim$ mass) or larger $-\alpha$ ($\sim$ phase velocity).
Roughly speaking, the natural confinement within one stable state tends to be realized in massless and/or light particles with sufficiently high traveling speeds as high as the speed of light.
The detailed mechanism of both enhancement and suppression of the transition due to the wave effect is discussed in our forthcoming paper.

\bibliographystyle{amsplain}

\end{document}